\documentclass[final]{siamltex}

\usepackage{amsfonts,amsmath,amssymb}
\usepackage{mathrsfs,mathtools}
\usepackage{enumerate}
\usepackage{hyperref}
\usepackage{esint}
\usepackage{graphicx}
\usepackage{bm}
\usepackage{cite}
\usepackage{xcolor}

\DeclareMathAlphabet{\mathpzc}{OT1}{pzc}{m}{it}

\title{A note on models for anomalous phase-change processes\thanks{This work has been partially supported by the Consejo Nacional de Investigaciones Cient\'ificas y T\'ecnicas (CONICET) of Argentina within a Postdoctoral Grant at Universidad Austral, Rosario, Argentina.
}}
\author{Andrea N. Ceretani \thanks{IMAS-UBA-CONICET, Intendente Guiraldes 2160, Capital Federal 1428, Argentina, \texttt{aceretani@dm.uba.ar}, and Escuela de Ciencia y Tecnolog\'ia, Universidad Nacional de San Mart\'in, Mart\'in de Irigoyen 3100, San Mart\'in 1650, Argentina, \texttt{aceretani@unsam.edu.ar}}}

\begin{document}  

\maketitle

\begin{abstract}
We review some fractional free boundary problems that were recently considered for modeling anomalous phase-transitions. All problems are of Stefan type and involve fractional derivatives in time according to Caputo's definition. We survey the assumptions from which they are obtained and observe that the problems are nonequivalent though all of them reduce to a classical Stefan problem when the order of the fractional derivatives is replaced by one. We further show that a simple heuristic approach built upon a fractional version of the energy balance and the classical Fourier's law leads to a natural generalization of the classical Stefan problem in which time derivatives are replaced by fractional ones.  

\end{abstract}



\section{Introduction}\label{Section:Intro}
This note addresses time-fractional free boundary problems of Stefan type that were recently considered in the literature as possible models for anomalous phase-transitions. Some of them consist simply of a Stefan problem in which the usual time derivatives are replaced by fractional ones, whereas others involve extra assumptions or include new terms when compared with a classical Stefan problem. All of them reduce to a usual Stefan problem when the fractional order of time derivatives is replaced by one.

Our interest is on phase-transitions that exhibit a characteristic diffusion time proportional to $t^{\gamma/2}$ with $\gamma\in(0,1)$, instead of the classical behavior determined by $t^{1/2}$. The anomalous diffusion is presumed to be caused by the effects of the process past history on its present state, and related free boundary problems are investigated.

The question that motivates this note is: Is there any fractional version of the classical Stefan problem that is more suitable in applications than the others? In order to answer this question, if possible, we investigate the physical hypothesis or mathematical assumptions leading to the fractional problems addressed here.

We consider a melting process for a phase-change material occupying a physical region $\Omega$ that is initially insulated. For simplicity, we assume the material to be initially solid at the phase change temperature, which we suppose equal to zero. 
We assume that the phase-change process is originated by a constant temperature prescribed at some part of the domain boundary, and that a liquid region is separated from a solid one by a sharp interface at any time. The melting process is described by the evolution in time of the temperature distribution and the location of the melt front.

To present ideas in a simple way, we further assume that $\Omega$ is a thin long cylinder with a heated base, hence the phase-change problem can be formulated in the one-dimensional semi-infinite domain $(0,\infty)$.

Consider that the lateral boundary of $\Omega$ remains insulated for all times. Then, the simplest model for the melting process described above is given by a {\em Stefan problem} \cite{AlSo1993}. Let $s:(0,\infty)\to\mathbb{R}$ be the location of the interface (free boundary) and $u:(0,\infty)\to\mathbb{R}$ be the material temperature. We denote by $u_\ell$ the temperature in the liquid region, thus 
\begin{equation}\label{Temperature}
u(x,t)=u_\ell(x,t)\,\,\, \mathrm{if}\,\,\, 0<x<s(t)\quad\text{and}\quad
u(x,t)=0\,\,\,\mathrm{otherwise},\,\,\,\text{for all }t>0.
\end{equation}
The Stefan problem establishes that the liquid temperature $u_\ell$ is described by the 
\begin{equation}\label{HeatEq}
\text{\em Heat equation}:\qquad\frac{\partial u_\ell}{\partial t}(x,t)=\alpha\frac{\partial^2 u_\ell}{\partial x^2}(x,t)\qquad 0<x<s(t),\quad t>0,
\end{equation}  
and the interface movement is characterized by the 
\begin{equation}\label{StefanCond}
\text{\em Stefan condition}:\qquad\rho l\dot{s}(t)=-k\frac{\partial u_\ell}{\partial x}(s(t)^-,t)\qquad t>0.
\end{equation}
The coefficient $\alpha$ in \eqref{HeatEq} is the thermal diffusivity, given by $\alpha:=\frac{k}{\rho c}$. The parameters $\rho$, $k$, $c$, and $l$ represent the material density, the thermal conductivity, the specific heat, and the latent heat in the liquid phase. All of them are assumed to be constant. The remaining conditions in the Stefan problem are
\begin{align}\label{Conditions}
s(0)=0,\quad 
u_\ell(x,0)=0\quad x>0,\quad u_\ell(0,t)=u_0\quad\text{and}\quad u_\ell(s(t),t)=0\quad t>0,
\end{align}
where $u_0>0$ is a constant temperature prescribed at the fixed boundary $x=0$.

The heat equation \eqref{HeatEq} and the Stefan condition \eqref{StefanCond} are derived from the following physical assumptions:
\begin{align}
\label{FirstPrThermod}&\text{\em First principle of thermodynamics}:\notag\\[0.25cm]
&\frac{\partial }{\partial t}\int_a^be(x,t)\,dx=\,q(a,t)-q(b,t)\qquad\text{for all $t>0$ and every $(a,b)\subset(0,\infty)$},\\[0.25cm]
\label{FourierLaw}&\text{\em Fourier's law}:\notag\\[0.25cm]
&q(x,t)=\,-k\frac{\partial u_\ell}{\partial x}(x,t)\qquad \text{for all\quad$0<x<s(t)$\quad and\quad$t>0$},
\end{align}
where $e$ is the energy per unit mass and $q$ is the heat flux. The latter is assumed to be zero in the solid phase.

Assumption \eqref{FirstPrThermod} implies:
\begin{align}
\label{EnergyConservLaw}\text{\em Energy conservation law}:&
\quad\,\,\,\frac{\partial e}{\partial t}(x,t)=\,\frac{\partial q}{\partial x}(x,t)\qquad 0<x<s(t),&t>0,&\\[0.25cm]
\label{GlobalEnergyLaw}\text{\em Global energy conservation}:&\quad\,\,\,
\frac{\partial}{\partial t}\int_0^\infty e(x,t)\,dx=\,q(0,t)\quad&t>0.&
\end{align}
Notice that $q(\infty,t)=0$ for all $t>0$. Considering that $e=\rho c u_\ell+\rho l$ in the liquid region, equation \eqref{HeatEq} arises from a plain combination of \eqref{FourierLaw} and \eqref{EnergyConservLaw}. Condition \eqref{StefanCond} is then obtained from the global balance \eqref{GlobalEnergyLaw}. A detailed derivation of the Stefan problem can be found in \cite{AlSo1993}.

It is worth mentioning that liquid infiltration in unsaturated porous media can be analogously described by a Stefan-like problem in which condition \eqref{StefanCond} characterizes the movement of the wet front, see e.g. \cite{EsFaMi2000,GiMa2001}.

The Stefan problem predicts a melt front advancement given by $s(t)\sim t^{1/2}$. However, interface movement given by $s(t)\sim t^{\gamma/2}$ with $\gamma\neq 1$ has been observed in phase-transitions and liquid infiltration in non-homogeneous domains, \cite{KuLa2001,FiHiLoVo2016,Vo2015,Vo2016}. For example, numerical experiments by Voller \cite{Vo2016} (see also \cite{Vo2015}) show phase-change processes with an interface location given by $s(t)\sim t^{\gamma/2}$ with $\gamma\in(0,1)$ in two dimensional domains in which the characteristic lenght of the heterogeneities is comparable with the characteristic lenght of the domain itself. Examples of diffusion-like processes with characteristic time of diffusion proportional to $t^{\gamma/2}$ with $\gamma\neq 1$ abound in real world applications, see e.g. the monograph \cite{MeKl2000} and references therein.

At present, there exists a vast literature on modeling anomalous transport processes in fixed domains by fractional differential equations. A recent survey is given in \cite{Va2017}. Equations involving fractional integrals or derivatives of non-integer order in time have shown to be suitable models for anomalous transport with diffusion time scale $\sim t^{\gamma/2}$, $\gamma\in(0,1)$, and a comprehensive mathematical theory for them is nowadays available, see e.g. \cite{Po1999,KiSrTr2006}. However, modeling anomalous phase-transitions by fractional free boundary problems is, so far, less developed. For example, several (nonequivalent) time-fractional Stefan-like problems were proposed as possible models for one-dimensional phase-change processes that exhibit a characteristic time of diffusion $\sim t^{\gamma/2}$ with $\gamma\in(0,1)$, e.g. \cite{VoFaGa2013,LiXu2009,RoBoTa2018}.

Fractional models usually arise as a consequence of the physical assumptions that, hopefully, describe the anomalous transport process under consideration. In section \ref{Section:Survey} we make a brief review of time-fractional Stefan-like problems recently studied in the literature, survey the assumptions from which they are obtained, and observe that they are nonequivalent. In particular, the latter suggests that the way in which ``memory'' is accounted for in a model is a sensitive issue.

Furthermore, in section \ref{Section:FractionalFBP} we derive a time-fractional Stefan-like problem from a heuristic approach based on a time-fractional version of condition \eqref{FirstPrThermod} and the classical Fourier's law. 
We obtain a fractional ``model" analogous to the classical Stefan problem since we recover the same governing equation \eqref{HeatEq} and condition on the interface \eqref{StefanCond} now with time-fractional derivatives instead of classical ones. This approach is motivated by the observation that the classical continuity equation may no longer be suitable in highly heterogeneous domains, see \cite{WhMe2008}, where transport processes may experience ``waiting times" do to particle trapping, see \cite{Vo2015}. Time-fractional conservation equations were also considered in \cite{CoMe1997}, where the relation with non-local transport theory with memory effects is discussed. In this way, we provide a new approach to derive some fractional problems that are of interest in pure and applied fields, \cite{LiXuWa2008,LiXuJi2009,DaRa2010,Vo2010,DaKuKu2011,SiKuRa2011,
At2012,RaSiKu2013,BlKl2015,LiWaZh2013,VoMiDa2016,Vo2014,
KuRiRy2017}.

\section{Brief review on time-fractional Stefan problems}\label{Section:Survey}
At present, different approaches at dealing with anomalous phase-transitions by time-fractional Stefan-like problems coexist in the literature. This section is devoted to a brief overview of them and provides a shortcut on the available research on the subject.

The fractional problems discussed below may seem, on the surface, pretty similar. However, the assumptions from which they are obtained are quite different from the physical point of view. In particular, this leads to nonequivalent models, as we show here.

Let $\gamma\in(0,1)$. We denote by $\frac{\partial^\gamma}{\partial t^\gamma}$ the left-sided Caputo fractional derivative of order $\gamma$, given by
\begin{equation}\label{Caputo}
\frac{\partial^\gamma f}{\partial t^\gamma}(t)=\frac{1}{\Gamma(1-\gamma)}\int_0^t(t-t')^{-\gamma}\frac{df}{dt'}(t')\,dt' \qquad t\geq 0,
\end{equation} 
for any $f\in C^1([0,\infty))$, where $\Gamma$ is the Gamma function. 
We notice that formula \eqref{Caputo} can be naturally extended to the larger class of absolutly continuous functions and refer to \cite{KiSrTr2006} for a detailed exposition on the Caputo fractional derivative. 

In \cite{VoFaGa2013}, Voller, Falcini, and Garra proposed a fractional model for the melting process described in section \ref{Section:Intro}, on a general bounded domain $\Omega$. They consider that $\Omega$ is initially insulated, and allow it to loose heat through the liquid boundary once the melting process has begun. The model is derived from a global energy balance and a constitutive equation that accounts for memory by means of a nonlocal-in-time definition of the heat flux.

The argument developed in \cite{VoFaGa2013} still holds true when the spatial domain is $(0,\infty)$, as considered in section \ref{Section:Intro}. For simplicity, we describe the ideas in \cite{VoFaGa2013} for this special case.

The model in \cite{VoFaGa2013} is derived from the global energy balance
\begin{equation}\label{VoFaGa-GobalBalanceEq}
\int_0^{s(t)} (\rho c u_\ell(x,t)+\rho l)\,dx=\int_0^t\hat{q}(0,t')\,dt'\qquad\forall\,t>0,
\end{equation}
where the (nonlocal) heat flux $\hat{q}$ is assumed to account for memory according to 
\begin{equation}\label{VoFaGa-ConstitutiveEq}
\hat{q}(x,t)=\frac{\tau^{1-\gamma}}{\Gamma(\gamma)}\frac{\partial}{\partial t}\int_0^t(t-t')^{\gamma-1}q(x,t')\,dt'\qquad 0<x<s(t),\quad t>0.
\end{equation}
In \eqref{VoFaGa-ConstitutiveEq}, $q$ denotes a (local) flux given by Fourier's law \eqref{FourierLaw} in the liquid phase and by zero in the solid one, and the parameter $\tau>0$ represents a constant relaxation time. Notice that assumption \eqref{VoFaGa-GobalBalanceEq} is weaker than \eqref{FirstPrThermod}, and that \eqref{VoFaGa-GobalBalanceEq} does not (necessarily) imply the local energy conservation law \eqref{EnergyConservLaw}.

The governing equation for the material temperature in the liquid phase and the condition describing the melt front movement obtained in \cite{VoFaGa2013} are given by
\begin{equation}\label{FracHeatEq}
\frac{\partial^\gamma u}{\partial t^\gamma}(x,t)=\alpha\tau^{1-\gamma}\frac{\partial^2 u_\ell}{\partial x^2}(x,t)\qquad 0<x<s(t),\quad t>0,
\end{equation}
and
\begin{equation}\label{FracStefanCond}
\begin{split}
\rho l\frac{\partial^\gamma s}{\partial t^\gamma}(t)=-\tau^{1-\gamma}k\frac{\partial u_\ell}{\partial x}(s(t)^-,t)\qquad t>0,
\end{split}
\end{equation}
respectively.

Differently from the usual approach given by \eqref{HeatEq} and \eqref{StefanCond}, the fractional model given by \eqref{FracHeatEq} and \eqref{FracStefanCond} prescribes a material temperature and a location of the interface coupled through the fractional differential operator in the governing equation. This becomes evident when we rewrite \eqref{FracHeatEq} taking into account the definition of $u$ given in \eqref{Temperature}:
\begin{equation}\label{FracHeatEq-withs}
\frac{\rho c}{\Gamma(1-\gamma)}\int\limits_{s^{-1}(x)}^t
(t-t')^{-\gamma}\frac{\partial u_\ell}{\partial t'}(x,t')\,dt'=k\tau^{\gamma-1}\frac{\partial^2 u_\ell}{\partial x^2}(x,t)\qquad 0<x<s(t),\quad t>0.
\end{equation}
In particular, this coupling breaks down the classical Neumann's method to find explicit similarity solutions to fractional Stefan-like problems involving \eqref{FracHeatEq} and \eqref{FracStefanCond}. To the best of our knowledge, the only explicit solution to this sort of problems is given in \cite{Vo2010} for the limit case of vanishing specific heat and predicts a melt font advancement proportional to $t^{\gamma/2}$. The aforementioned coupling also becomes challenging the weak solvability of free boundary problems involving \eqref{FracHeatEq} and \eqref{FracStefanCond}. Some advances in this direction were recently given in \cite{KuRiRy2017}. Looking for approximate solutions (both numerical and analytical) is an active area of research, \cite{LiWaZh2013,LiXuJi2009,DaKuKu2011,DaRa2010,SiKuRa2011,
RaSiKu2013,VoMiDa2016,BlKl2015,Vo2014}.

Free boundary problems involving fractional heat equation \eqref{FracHeatEq} and fractional Stefan condition \eqref{FracStefanCond} were considered in several works at dealing with anomalous diffusive-like processes, e.g. \cite{LiXuWa2008,LiXuJi2009,DaRa2010,Vo2010,DaKuKu2011,SiKuRa2011,
At2012,RaSiKu2013,BlKl2015,LiWaZh2013,VoMiDa2016,Vo2014}. A remarkable contribution from Voller, Falcini, and Garra was to identify a definition for the heat flux and a suitable classical global balance leading to a transparent generalization of the classical model to a fractional one. However, the presence of the time derivative in the right-hand side of the constitutive equation \eqref{VoFaGa-ConstitutiveEq} is still challenging to endow \eqref{VoFaGa-ConstitutiveEq} a physical meaning, see \cite{GuPi1968}. In next section \ref{Section:FractionalFBP} we shall see that equation \eqref{FracHeatEq} and condition \eqref{FracStefanCond} can also be obtained from a different approach based on a fractional version of the conservation principle \eqref{FirstPrThermod} and the classical Fourier's law \eqref{FourierLaw}.

Other approach to generalize Stefan problems to the fractional framework consists of replacing the time derivatives in \eqref{HeatEq} and \eqref{StefanCond} by fractional ones by means of the Caputo operator, {\em and simultaneously neglecting definition \eqref{Temperature} for the material temperature}. Stefan condition is still replaced by \eqref{FracStefanCond} and the heat equation \eqref{HeatEq} is now substituted by
\begin{equation}\label{FracHeatEq-bis}
\begin{split}
&\frac{\partial^\gamma u}{\partial t^\gamma}=\alpha\tau^{1-\gamma}\frac{\partial^2 u_\ell}{\partial t^{\gamma}}\quad 0<x<s(t),\quad t>0,\\[0.25cm]
&\text{where\quad $u(x,t)=u_\ell(x,t)$\quad for\quad $0<x<s(t),\quad t>0$},\\[0.25cm]
&\text{(no condition on $u(x,t)$ for $x>s(t)$)}.
\end{split}
\end{equation}

From the physical point of view, removing the assumption $u(x,t)=0$ for $x>s(t)$, $t>0$, completely changes the problem nature. However, satisfactory results were obtained following this approach, see e.g. \cite{LiXu2004,Vo2014,DaKuKu2011,DaRa2010,RaKu2013,
LiSu2015,LiXuWa2007} (see also \cite{LiXu2009,RoTa2014,Ta2015-b,RoSa2014,RoSa2013,Ro2016,
Ro2017,CeTa2017,RoTa2018}). 

An attractive feature of this second type of fractional problems is that they admit explicit solutions, \cite{LiXu2009,RoSa2013}.
While on the surface \eqref{FracHeatEq} and \eqref{FracHeatEq-bis} seem similar, the second equation does not (necessarily) couple $u$ and $s$ through the differential operator in the left-hand side. This feature is key to find explicit solutions by Neumann's method, as pointed out by Liu and Xu in \cite{LiXu2009} and by Roscani and Santillan Marcus in \cite{RoSa2013}. Below, we briefly describe the argument in \cite{LiXu2009,RoSa2013} for the problem given by \eqref{FracHeatEq-bis}, \eqref{FracStefanCond}, and \eqref{Conditions}.

The first step is to notice that the differential equation in \eqref{FracHeatEq-bis} can be written as 
\begin{equation}\label{FractHeatEq-bis-bis}
\frac{\partial^\gamma u}{\partial t^\gamma}(x,t)=\alpha\tau^{1-\gamma}\frac{\partial^2 u}{\partial x^2}(x,t)\qquad x>0,\quad t>0,
\end{equation} 
which admits solutions $u$ of the form
\begin{equation}\label{uSolution-bis}
u(x,t)=a+bW\left(-\frac{x}{\lambda t^{\gamma/2}},-\frac{\gamma}{2},1\right)
\qquad x>0,\quad t>0,\quad(a, b\in \mathbb{R})
\end{equation}
where $\lambda:=\tau^{(1-\gamma)/2}\alpha^{\gamma/2}$ and $W$ is the Wright function, given by
\begin{equation*}
W(z,\mu,\nu)=\displaystyle\sum_{k=0}^{\infty}\frac{z^k}{k!\,\Gamma(\mu k+\nu)}
\qquad z\in\mathbb{R},\quad \mu>-1,\quad\nu\in\mathbb{R}.
\end{equation*} 
Then, a solution of similarity type to \eqref{FracHeatEq-bis}-\eqref{FracStefanCond}-\eqref{Conditions} is obtained by considering $u$ given by \eqref{uSolution-bis}, $s$ in the form $s(t)=\sigma\lambda t^{\gamma/2}$ ($\sigma\in\mathbb{R}$),
and selecting $a$, $b$, and $\sigma$ in order to fulfill \eqref{FracStefanCond} and \eqref{Conditions}. In this way, the following solution is obtained:
\begin{align}
\label{u}&u_\ell(x,t)=u_0\left(1-\frac{1-W\left(-\frac{x}{\lambda t^{\gamma/2}},-\frac{\gamma}{2},1\right)}{1-W\left(-\sigma,-\frac{\gamma}{2},1\right)}\right)&
&0<x<s(t),& t>0,&\\[0.25cm]
\label{s}
&s(t)=\sigma\lambda t^{\gamma/2}&
&&t>0,&
\end{align}
where the parameter $\sigma>0$ is defined as the solution to a specific transcendental equation (see details in \cite{LiXu2009,RoSa2013}). So far, the only explicit solutions to fractional free boundary problems involving \eqref{FracHeatEq-bis} and \eqref{FracStefanCond} are obtained from this method.

It is worth mentioning that the function $u$ given by \eqref{u} does not satisfy \eqref{FracHeatEq-withs}: Taking the time derivative of $u$ in \eqref{u}, we obtain
\begin{equation*}
\frac{\partial  u}{\partial t}(x,t)=\left(\frac{\gamma u_0 }{2\lambda}\right) \frac{W'\left(-\frac{x}{\lambda t^{\gamma/2}},-\frac{\gamma}{2},1\right)\,xt^{-1+\gamma/2}}{1-W(-\sigma,-\gamma/2,1)}<0\quad\text{for all}\quad x>0,\quad t>0,
\end{equation*} 
since $0<W(-x,-\gamma/2,1)<1$ for all $x> 0$ and $W(-x,-\gamma/2,1)$ is strictly decreassing for $x>0$ (see \cite{RoSa2013}). Here, $W'$ denotes the derivative of $W$ with respect to the first argument. In addition, we observe that $t'<t$ provided that $0<t'<s^{-1}(x)$, for any $0<x<s(t)$ and $t>0$. Therefore, 
\begin{equation*}
(t-t')^{-\gamma}\frac{\partial u}{\partial t'}(x,t')<0\qquad \text{for all}\quad 0<t'<s^{-1}(x),\quad 0<x<s(t),\quad t>0.
\end{equation*}
From this we observe that $u$, given by \eqref{u}, does not satisfy \eqref{FracHeatEq-withs} since   
\begin{equation*}
\begin{split}
&\frac{\rho c}{\Gamma(1-\gamma)}\int\limits_{s^{-1}(x)}^t
(t-t')^{-\gamma}\frac{\partial u_\ell}{\partial t'}(x,t')\,dt'
\\[0.25cm]
&\quad=\tau^{1-\gamma}k\frac{\partial^2 u_\ell}{\partial x^2}(x,t)-
\frac{\rho c}{\Gamma(1-\gamma)}\int\limits_0^{s^{-1}(x)}
(t-t')^{-\gamma}\frac{\partial u}{\partial t'}(x,t')\,dt',
\end{split}
\end{equation*}
for all $0<x<s(t)$ and $t>0$ (see \eqref{FracHeatEq-bis}), and  the second term in the right-hand side does not vanish identically.

In particular, the above shows that the fractional Stefan problems given by \eqref{FracHeatEq-bis}-\eqref{FracStefanCond}-\eqref{Conditions} with $u$ given by \eqref{Temperature}, and \eqref{FracHeatEq}-\eqref{FracStefanCond}-\eqref{Conditions} are not equivalent for $\gamma\in(0,1)$. In addition, we observe that solutions to these problems may account for memory in a different way since the latter neglects \eqref{Temperature}. Further relations between the two problems seem, so far, not to have been considered in the literature (e.g., how ``close" are they? does it depend on the value of $\gamma$?).

At present, most models at dealing with anomalous diffusion in domains with free or moving boundaries are like the (nonequivalent) problems \eqref{FracHeatEq}-\eqref{FracStefanCond}-\eqref{Conditions} with $u$ given by \eqref{Temperature}, or \eqref{FracHeatEq-bis}-\eqref{FracStefanCond}-\eqref{Conditions}.

Recently, Roscani, Bollati, and Tarzia introduced a third approach in \cite{RoBoTa2018} that leads to a new fractional Stefan-like problem for one-dimensional phase-transitions. The model in \cite{RoBoTa2018} is given by \eqref{Conditions} and
\begin{align}
\label{FracEq}&\frac{\partial^\gamma u}{\partial t^\gamma}(x,t)+
\frac{lc^{-1}(t-s^{-1}(x))^{-\gamma}}{\Gamma(1-\gamma)}=\tau^{1-\gamma}\alpha\frac{\partial^2 u_\ell}{\partial x^2}(x,t)&0<x<s(t),\quad t>0,\\[0.25cm]
\label{FracCond}&\rho ls'(t)=-\tau^{1-\gamma}k\frac{\hat{\partial}^{1-\gamma}}{\partial t^{1-\gamma}}u(s(t)^-,t)&t>0,
\end{align} 
where $u$ is defined by \eqref{Temperature} and $\frac{\hat{\partial}^{1-\gamma}}{\partial t^{1-\gamma}}$ denotes the Riemann-Liouville fractional derivative of order $1-\gamma$. We recall that the Riemann-Liouville fractional operator $\frac{\hat{\partial}^\delta}{\partial t^\delta}$, $\delta\in(-1,1)$, is defined by
\begin{equation}\label{Riemann-Liouville}
\frac{\hat{\partial}^\delta f}{\partial t^\delta}:=\left\{
\begin{array}{lcc}
\frac{1}{\Gamma(1-\delta)}\frac{\partial}{\partial t}\int_0^t(t-t')^{-\delta}f(t')\,d t'&\mathrm{if}&0\leq \delta<1,\\[0.25cm]
\frac{1}{\Gamma(-\delta)}\int_0^t(t-t')^{-\delta-1}f(t')\,d t'&\mathrm{if}&-1<\delta<0,
\end{array}
\right.
\end{equation}
see \cite{KiSrTr2006}. Observe that $\frac{\hat{\partial}^\delta}{\partial t^\delta}$ is a differential operator for $\delta\in(0,1)$, whereas it is an integral operator for $\delta\in(-1,0)$.

The fractional equation \eqref{FracEq} and condition \eqref{FracCond} are obtained from the non-classical constitutive equation \eqref{VoFaGa-ConstitutiveEq} in combination with a classical local energy balance. Notice that, even when the heat flux is assumed to account for memory in the same way in Roscani's et al. and Voller's et al. approaches, the former considers a local energy balance in the liquid region instead of the (weaker) global balance assumed in the second one, see \eqref{VoFaGa-GobalBalanceEq}. The existence of solutions to the problem introduced in \cite{RoBoTa2018} was not yet investigated.

\section{A heuristic approach leading to a time-fractional Stefan problem}\label{Section:FractionalFBP}
This section is aimed to provide another approach to derive the time-fractional Stefan problem obtained in \cite{VoFaGa2013} (see the first problem discussed in section \ref{Section:Survey}), which involves equation \eqref{FracHeatEq} and condition \eqref{FracStefanCond}. It consists of considering a fractional time derivative, instead of the usual one, in the classical approach described in section \ref{Section:Intro} (see \eqref{FirstPrThermod} and \eqref{FourierLaw}). 

All through this section, the temperature $u_\ell$ and the heat flux $q$ in the liquid phase are assumed to be sufficiently smooth to make the following computations meaningful. In addition, we suppose that $s$ is an increasing function (i.e., the melt front can only advance in time).

We consider the following mathematical generalization of the physical principle \eqref{FirstPrThermod}:
\begin{equation}\label{FracFirstPrThermd}
\frac{\partial^\gamma}{\partial t^\gamma}\int_a^be(x,t)\,dx=\,\tau^{1-\gamma}(q(a,t)-q(b,t))\qquad\text{for all $t>0$ and every $(a,b)\subset(0,\infty)$},
\end{equation}
where $e$ is the energy per unit mass and $q$ is the heat flux, which we suppose given by the Fourier's law \eqref{FourierLaw} in the liquid region and by zero in the solid one. The parameter $\tau>0$ in \eqref{FracFirstPrThermd} represents a constant relaxation time.

We assume that the energy is given by the enthalpy of the system, that is,
\begin{equation*}
e(x,t)=\rho c u_\ell(x,t)+\rho l\,\,\,\mathrm{if}\,\,\,0<x\leq s(t),\,t>0,\quad\text{and}\quad
e(x,t)=0\,\,\,\mathrm{otherwise}.
\end{equation*} 
We further note that \eqref{FracFirstPrThermd} implies the following fractional version of the energy conservation law \eqref{EnergyConservLaw}:
\begin{equation}\label{FracEnergyConservLaw}
\frac{\partial^\gamma e}{\partial t^\gamma}(x,t)=\,\tau^{1-\gamma}\frac{\partial q}{\partial x}(x,t)\qquad 0<x<s(t),\quad t>0,
\end{equation}
which combined with the Fourier's law \eqref{FourierLaw} yields 
\begin{equation}
\frac{\partial^\gamma u}{\partial t^\gamma}(x,t)=\alpha\tau^{1-\gamma}\frac{\partial^2 u_\ell}{\partial x^2}(x,t)\qquad 0<x<s(t),\quad t>0.
\end{equation}
Here, we have taken into account that $\frac{\partial^\gamma(\rho c)}{\partial t^\gamma}=0$ since $\rho c$ is constant.

Assumption \eqref{FracFirstPrThermd} also implies the following fractional version of the global energy balance \eqref{GlobalEnergyLaw}:
\begin{equation}\label{FracGlobalEnergyBalance}
\frac{\partial^\gamma}{\partial t^\gamma}\int_0^\infty e(x,t)\,dx=\tau^{1-\gamma}q(0,t)\qquad t>0.
\end{equation}
We observe:
\begin{equation*}
\begin{split}
&\frac{\partial^\gamma }{\partial t^\gamma}\int\limits_0^\infty e(x,t)\,dx=\,
\frac{\partial^\gamma }{\partial t^\gamma}\int\limits_0^{s(t)} e(x,t)\,dx\\[0.25cm]
&\quad=\frac{1}{\Gamma(1-\gamma)}\displaystyle\int_0^t
(t-t')^{-\gamma}\frac{\partial}{\partial t'}
\left(\int\limits_0^{s(t)}(\rho c u_\ell(x,t')+\rho l)\,dx\right)dt'\\[0.25cm]
&\quad=
\frac{\rho c}{\Gamma(1-\gamma)}\int\limits_0^t
(t-t')^{-\gamma}
\left(u_\ell(s(t')^-,t')\dot{s}(t')+
\int\limits_0^{s(t')}\frac{\partial u_\ell}{\partial t'}(x,t')\,dx\right)\,dt'+\rho l\frac{\partial^\gamma s}{\partial t^\gamma}(t)\\[0.25cm]
&\quad=
\frac{\rho c}{\Gamma(1-\gamma)}\int\limits_0^t
\int\limits_0^{s(t')}
(t-t')^{-\gamma}
\frac{\partial u_\ell}{\partial t'}(x,t')\,dx\,dt'+\rho l\frac{\partial^\gamma s}{\partial t^\gamma}(t)\\[0.25cm]
&\quad=\frac{\rho c}{\Gamma(1-\gamma)}\int\limits_0^{s(t)}
\int\limits_{s^{-1}(x)}^t
(t-t')^{-\gamma}\frac{\partial u_\ell}{\partial t'}(x,t')\,dt'\,dx+\rho l\frac{\partial^\gamma s}{\partial t^\gamma}(t)\\[0.25cm]
&\quad=\rho c\int\limits_0^{s(t)}\frac{\partial^\gamma u}{\partial t^\gamma}(x,t)\,dx+\rho l\frac{\partial^\gamma s}{\partial t^\gamma}(t).
\end{split}
\end{equation*}
We have used here that $s$ is invertible and that the interface is at the melting temperature, see \eqref{Conditions}. Then, taking into account that $u_\ell$ satisfies \eqref{FracHeatEq}, we have:
\begin{equation*}
\begin{split}
&\frac{\partial^\gamma}{\partial t^\gamma}\int\limits_0^\infty e(x,t)\,dx=\,
\tau^{1-\gamma}k\int\limits_0^{s(t)}\frac{\partial^2 u_\ell}{\partial x^2}(x,t)\,dx+\rho l\frac{\partial^\gamma s}{\partial t^\gamma}(t)\\[0.25cm]
&\quad=
\tau^{1-\gamma}k\frac{\partial u_\ell}{\partial x}(s(t)^-,t)-
\tau^{1-\gamma}k\frac{\partial u_l}{\partial x}(0^+,t)+\rho l\frac{\partial^\gamma s}{\partial t^\gamma}(t)\\[0.25cm]
&\quad=
\tau^{1-\gamma}k\frac{\partial u_\ell}{\partial x}(s(t)^-,t)+\tau^{1-\gamma}q(0^+,t)+\rho l\frac{\partial^\gamma s}{\partial t^\gamma}(t).
\end{split}
\end{equation*}
Therefore, \eqref{FracGlobalEnergyBalance} yields 
\begin{equation}
\begin{split}
\rho l\frac{\partial^\gamma s}{\partial t^\gamma}(t)=-\tau^{1-\gamma}k\frac{\partial u_\ell}{\partial x}(s(t)^-,t)\qquad t>0.
\end{split}
\end{equation}

Hence, the fractional generalization \eqref{FracFirstPrThermd} of the physical condition \eqref{FirstPrThermod} together with Fourier's law \eqref{FourierLaw} predicts a material temperature in the liquid region described by the fractional heat equation \eqref{FracHeatEq} and a melt front advancement characterized by the fractional Stefan condition \eqref{FracStefanCond}.

The argument presented here extends naturally to the case when the initial material temperature is less than the melting one, thus two-phase fractional Stefan-like problems can be similarly obtained (see \cite{AlSo1993} for the classical derivation). Models for problems with interfaces of non-zero thickness can be analogously generalized to the fractional context, (see \cite{SoWiAl1982,Ta1990} for classical derivations).

Finally, we notice that the assumptions in \cite{RoBoTa2018} (see the third problem discussed in section \ref{Section:Survey}) do not imply those considered here. In fact: The (nonlocal) heat flux \eqref{VoFaGa-ConstitutiveEq} in combination with a classical energy conservation law in the liquid phase yields
\begin{equation*}
\frac{\partial e}{\partial t}(x,t)=\frac{\partial \hat{q}}{\partial x}(x,t)\qquad 0<x<s(t),\quad t>0.
\end{equation*}
Applying the Riemann-Liouville integral $\frac{\partial^{\gamma-1}}{\partial t^{\gamma-1}}$ to both sides and taking into account the definition of $\hat{q}$ given in \eqref{VoFaGa-ConstitutiveEq}, we have
\begin{equation*}
\frac{\partial^\gamma e}{\partial t^\gamma}(x,t)=\frac{\hat{\partial}^{\gamma-1}}{\partial t^{\gamma-1}}\frac{\partial }{\partial x}\frac{\hat{\partial}^{1-\gamma}}{\partial t^{1-\gamma}}q(x,t)\qquad 0<x<s(t),\quad t>0,
\end{equation*}
which in general differs from \eqref{FracEnergyConservLaw} since, as pointed out in \cite{RoBoTa2018}, we cannot interchange the first integral operator with the classical spacial derivative in the right-hand side due to the loss of differentiability of the heat flux $\hat{q}$ along the interface.

\section{Concluding remarks}
We considered time-fractional free boundary problems of Stefan type that were studied in the recent literature, and discussed them in the light of modeling anomalous phase-transitions. All problems involve time-fractional differential operators of order $\gamma\in(0,1)$ and reduce to a classical Stefan problem when time derivatives are considered in the usual sense ($\gamma=1$). However, we saw that they may no longer be equivalent when $\gamma\in(0,1)$: Problems obtained by classical (local or global) conservation energy equations and nonlocal constitutive equations for the heat flux may not be equivalent among them, nor with problems derived from non-classical conservation laws in combination with local constitutive equations. This suggests that fractional models should be carefully chosen in order to fit the physical problem under consideration. We further presented a simple heuristic approach based on a fractional version of the conservation of energy and the classical Fourier's law that leads to a transparent generalization of the classical one-phase Stefan problem, in which the usual time derivatives are replaced by fractional ones according to Caputo's definition. This provides another interpretation of this kind of problems, in addition to the already known meaning in terms of nonlocal transport theory. The approach proposed here can be naturally extended to two-phases problems or phase-transitions with non-zero thickness interfaces.

\section*{Acknowledgements}
The author thanks Roberto Garra (La Sapienza, Rome, Italy), Federico Falcini (Consiglio Nazionale delle Ricerche, Rome, Italy), Vaughan Voller (University of Minnesota, Minneapolis, USA), Domingo Tarzia (Universidad Austral - CONICET, Rosario, Argentina), and Sabrina Roscani (Universidad Austral - CONICET, Rosario, Argentina) for fruitful discussions.

\bibliographystyle{plain} 
\bibliography{References}

\begin{thebibliography}{10}

\bibitem{AlSo1993}
V.~Alexiades and A.~D. Solomon.
\newblock {\em Mathematical modeling of melting and freezing processes}.
\newblock Hemisphere Publishing Corporation, Washington, 1993.

\bibitem{At2012}
C.~Atkinson.
\newblock Moving boundary problems for time fractional and composition
  dependent diffusion.
\newblock {\em Fractional Calculus and Applied Analysis}, 15(2):207--221, 2012.

\bibitem{BlKl2015}
M.~Blasik and M.~Klimek.
\newblock Numerical solution for the one phase 1{D} fractional {S}tefan problem
  using front fixing method.
\newblock {\em Mathematical Methods in the Applied Sciences}, 38:3214--3228,
  2015.

\bibitem{CeTa2017}
A.~N. Ceretani and D.~A. Tarzia.
\newblock Determination of two unknown thermal coefficients through and inverse
  one-phase fractional {S}tefan problems.
\newblock {\em Fractional Calculus and Applied Analysis}, 20(2):14--17, 2017.

\bibitem{CoMe1997}
A.~Compte and R.~Metzler.
\newblock The generalized {C}attaneo equation for the description of anomalous
  transport processes.
\newblock {\em Journal of Physics A: Mathematical and General}, 30:7277--7289,
  1997.

\bibitem{DaKuKu2011}
S.~Das, R.~Kumar, and P.~Kumar~Gupta.
\newblock Analytical approximate solution of space-time fractiona diffusion
  equation with a moving boundary condition.
\newblock {\em Zeitschrift f\"ur Naturforschung}, 66a:281--288, 2011.

\bibitem{DaRa2010}
S.~Das and Rajeev.
\newblock Solution of fractional diffusion equation with a moving boundary
  condition by variational iteration method and {A}domian decomposition method.
\newblock {\em Zeitschrift f\"ur Naturforschung}, 65a:793--799, 2010.

\bibitem{EsFaMi2000}
M.S. Espedal, A.~Fasano, and A.~Mikeli\'c.
\newblock {\em Filtration in porous media and industrial application: lectures
  given at the 4th session of the {C}entro {I}nternazionale {M}atematico
  {E}stivo}.
\newblock Springer, 2000.

\bibitem{FiHiLoVo2016}
N.~Filipovitch, K.~M. Hill, A.~Longjas, and V.~R. Voller.
\newblock Infiltration experiments demonstrate an explicit connection between
  heterogeneity and anomalous diffusion behavior.
\newblock {\em Water Resources Research}, 52(7):5167--5178, 2016.

\bibitem{GiMa2001}
R.~Gianni and P.~Mannucci.
\newblock A free boundary problem in an absorbing porous material with
  saturation dependent permeability.
\newblock {\em Nonlinear Differential Equantions and Applications}, 8:219--235,
  2001.

\bibitem{GuPi1968}
M.~E. Gurtin and A.~C. Pipkin.
\newblock A general theory of heat conduction with finite wave speeds.
\newblock {\em Archive for Rational Mechanics and Analysis}, 31:113--126, 1968.

\bibitem{KiSrTr2006}
A.~Kilbas, H.~Srivastava, and H.~Trujillo.
\newblock {\em Theory and applications of fractional differential equations}.
\newblock Elsevier, Amsterdam, 2006.

\bibitem{KuRiRy2017}
A.~Kubica, P.~Ribka, and K.~Ryszewska.
\newblock Weak solutions of fractional differential equations in non
  cylindrical domains.
\newblock {\em Nonlinear Analysis: Real World Applications}, 36:154--182, 2017.

\bibitem{KuLa2001}
M.~K\"untz and P.~Lavall\'ee.
\newblock Experimental evidence and theoretical analysis of anomalous diffusion
  during water infiltration in porous building materials.
\newblock {\em Journal of Physics D: Applied Physics}, 34:2547--2554, 2001.

\bibitem{LiSu2015}
X.~Li and X.~Sun.
\newblock Similarity solutions for phase-change problems with fractional
  governing equations.
\newblock {\em Applied Mathematics Letters}, 45:7--11, 2015.

\bibitem{LiWaZh2013}
X.~Li, S.~Wang, and M.~Zhao.
\newblock Two methods to solve a fractional single phase moving boundary
  problem.
\newblock {\em Central European Journal of Physics}, 11(10):1387--1391, 2013.

\bibitem{LiXuJi2009}
X.~Li, M.~Xu, and X.~Jiang.
\newblock Homotophy perturbation method to time-fractional diffusion equation
  with a moving boundary condition.
\newblock {\em Applied Mathematics and Computations}, 208:434--439, 2009.

\bibitem{LiXuWa2007}
X.~Li, M.~Xu, and S.~Wang.
\newblock Analytical solution to the moving boundary problems with
  space-time-fractional derivatives in drug release devises.
\newblock {\em Journal of Physics A: Mathematical and Theoretical},
  40:12131--12141, 2007.

\bibitem{LiXuWa2008}
X.~Li, M.~Xu, and S.~Wang.
\newblock Scale-invariant solutions to partial differential equations of
  fractional order with a moving boundary condition.
\newblock {\em Journal of Physics A: Mathematical and Theoretical}, 41:155202,
  2008.

\bibitem{LiXu2004}
X.~Liu and M.~Xu.
\newblock An exact solution to the moving boundary problem with fractional
  anomalous diffusion in drug release devises.
\newblock {\em ZAMM - Journal of Applied Mathematics and Mechanics},
  84(1):22--28, 2004.

\bibitem{LiXu2009}
X.~Liu and M.~Xu.
\newblock Some exact solutions to {S}tefan problems with fractional
  differential equations.
\newblock {\em Journal of Mathematical Analysis and Applications},
  351:536--542, 2009.

\bibitem{MeKl2000}
R.~Metzler and J.~Klafter.
\newblock The random walk's guide to anomalous diffusion: a fractional dynamics
  approach.
\newblock {\em Physics Reports}, 339:1--77, 2000.

\bibitem{Po1999}
I.~Podlubny.
\newblock {\em Fractional Differential Equations}.
\newblock Academic Press, San Diego, 1971.

\bibitem{RaKu2013}
Rajeev and M.S. Kushwaha.
\newblock Homotopy perturbation method for a limit case problem governed by
  fractional diffusion equation.
\newblock {\em Applied Mathematical Modelling}, 37:3589--3599, 2013.

\bibitem{RaSiKu2013}
Rajeev, M.~Singh~Kushwaha, and A.~Kumar.
\newblock An approximate solution to a moving boundary problem with space-time
  fractional derivative in fluvio-deltaic sedimentation process.
\newblock {\em Ain Shams Engineering Journal}, 4:889--895, 2013.

\bibitem{Ro2016}
S.~Roscani.
\newblock Hopf lemma for the fractional diffusion operator and its application
  to a fractional free-boundary problem.
\newblock {\em Journal of Mathematical Analysis and Applications},
  434:125--135, 2016.

\bibitem{Ro2017}
S.~Roscani.
\newblock Moving-boundary problems for the time-fractional diffusion equation.
\newblock {\em Electronic Journal of Differential Equations}, 2017(44):1--12,
  2017.

\bibitem{RoBoTa2018}
S.~Roscani, J.~Bollati, and D.A. Tarzia.
\newblock A new mathematical formulation for a phase change problem with a
  memory flux.
\newblock {\em Chaos, Solitons and Fractals}, 116:340--347, 2018.

\bibitem{RoSa2013}
S.~Roscani and E.~Santillan~Marcus.
\newblock Two equivalent {S}tefan's problems for the time-fractional-diffusion
  equation.
\newblock {\em Fractional Calculus in Applied Analysis}, 16(4):802--815, 2013.

\bibitem{RoSa2014}
S.~Roscani and E.~Santillan~Marcus.
\newblock A new equivalence of {S}tefan's problems for the
  time-fractional-diffusion equation.
\newblock {\em Fractional Calculus in Applied Analysis}, 17(2):371--381, 2014.

\bibitem{RoTa2014}
S.~Roscani and D.~A. Tarzia.
\newblock A generalized {N}eumann solution for the two-phase fractional
  {L}am\'e-{C}lapeyron-{S}tefan problem.
\newblock {\em Advances in Mathematical Sciences and Applications},
  24(2):237--249, 2014.

\bibitem{RoTa2018}
S.~Roscani and D.~A. Tarzia.
\newblock Explicit solution for a two-phase fractional {S}tefan problem with a
  heat flux condition at the fixed face.
\newblock {\em Computational and Applied Mathematics}, 37(5):4757--4771, 2018.

\bibitem{SiKuRa2011}
J.~Singh, P.~K. Gupta, and K.~N. Rai.
\newblock Variational iteration method to solve moving boundary problem with
  temperature dependent physical properties.
\newblock {\em Thermal Sciences}, 15, Suppl. 2:S229--S239, 2011.

\bibitem{SoWiAl1982}
A.~D. Solomon, D.~G. Wilson, and V.~Alexiades.
\newblock A mushy zone model with an exact solution.
\newblock {\em Letters in Heat and Mass Transfer}, 9:319--324, 1982.

\bibitem{Ta1990}
D.~A. Tarzia.
\newblock Neumann-like solution for the two-phase {S}tefan problem with a
  simple mushy zone model.
\newblock {\em Computational and Applied Mathematics}, 9-3:201--211, 1990.

\bibitem{Ta2015-b}
D.~A. Tarzia.
\newblock Determination of one unknown thermal coefficient through the
  one-phase fractional {L}am\'e-{C}lapeyron-{S}tefan problem.
\newblock {\em Applied Mathematics}, 6:2128--2191, 2015.

\bibitem{Va2017}
J.~L. V{\'{a}}zquez.
\newblock The mathematical theories of diffusion: Nonlinear and fractional
  diffusion.
\newblock In {\em Nonlocal and Nonlinear Diffusions and Interactions: New
  Methods and Directions}, pages 205--278. Springer International Publishing,
  2017.

\bibitem{VoMiDa2016}
C.~Vogl, M.~Miksis, and S.~Davis.
\newblock Moving boundary problems governed by anomalous diffusion.
\newblock {\em Proceedings of the Royal Society A}, 468:3348--3369, 2016.

\bibitem{Vo2010}
V.~R. Voller.
\newblock An exact solution of a limit case {S}tefan problem governed by a
  fractional diffusion equation.
\newblock {\em International Journal of Heat and Mass Transfer}, 53:5622--5625,
  2010.

\bibitem{Vo2014}
V.~R. Voller.
\newblock Fractional {S}tefan problems.
\newblock {\em International Journal of Heat and Mass Transfer}, 74:269--277,
  2014.

\bibitem{Vo2015}
V.~R. Voller.
\newblock A direct simulation demonstrating the role of spacial heterogeneity
  in determining anomalous diffusive transport.
\newblock {\em Water Resources Research}, 51(4):2119--2127, 2015.

\bibitem{Vo2016}
V.~R. Voller.
\newblock Computations of anomalous phase change.
\newblock {\em International Journal of Numerical Methods for Heat \& Fluid
  Flow}, 26(3/4):624--638, 2016.

\bibitem{VoFaGa2013}
V.~R. Voller, F.~Falcini, and R.~Garra.
\newblock Fractional {S}tefan problems exhibiting lumped and distributional
  latent-heat memory effects.
\newblock {\em Physical Review E}, 87:042401, 2013.

\bibitem{WhMe2008}
S.~W. Wheatcraft and M.~M. Meerschaert.
\newblock Fractional conservation of mass.
\newblock {\em Advances in Water Resources}, 31:1377--1381, 2008.

\end{thebibliography}
\newpage
\end{document}